\documentclass[a4paper]{article}
\usepackage{latexsym}
\usepackage{amsmath,amssymb,amsfonts}

\newtheorem{prethm}{{\bf Theorem}}

\newenvironment{thm}{\begin{prethm}{\hspace{-0.5
               em}{\bf.}}}{\end{prethm}}

\newtheorem{prepro}[prethm]{Proposition}

\newtheorem{prelem}[prethm]{Lemma}
\newenvironment{lem}{\begin{prelem}{\hspace{-0.5
               em}{\bf.}}}{\end{prelem}}

\newtheorem{precor}[prethm]{Corollary}
\newenvironment{cor}{\begin{precor}{\hspace{-0.5
               em}{\bf.}}}{\end{precor}}

\newtheorem{preremark}{{\bf Remark}}

\newtheorem{preexample}{{\bf Example}}

\newtheorem{preproof}{{\bf Proof.}}

\newenvironment{proof}[1]{\begin{preproof}{\rm
               #1}\hfill{$\Box$}}{\end{preproof}}

\title{\bf\ Distance-balanced closure of some graphs}
\input{epsf}
\author{N. Ghareghani $^{{\,\rm b, d}, \, 1, 2}$,  B. Manoochehrian $^{\,\rm c}$, M. Mohammad-Noori$^{\,\rm a,d}$\\
{\footnotesize {\em $^{\rm a}$Department of Mathematics, Statistics and Computer Science, University of Tehran}}\\
{\footnotesize {\em P.O. Box 14155-6455, Tehran, Iran}}\\
{\footnotesize{$^{ \textrm{b}}$\it Department of Mathematical Sciences, K.N. Toosi University of Technology}}\\
{\footnotesize{$^{ \textrm{c}}$\it Academic Center for Education, Culture and Research (ACECR), Tehran Branch, Tehran, Iran}}\\
{\footnotesize {\em $^{\rm d}$School of Mathematics, Institute for Research in Fundamental Sciences (IPM),}}\\
 {\footnotesize {\em P.O. Box 19395-5746, Tehran, Iran}}
 \\{\footnotesize Emails: ghareghani@ipm.ir, behzad@khayam.ut.ac.ir},
 \\{\footnotesize  morteza@ipm.ir, mnoori@khayam.ut.ac.ir}
 }
\begin{document}
\maketitle
\footnotetext[1]{\tt Corresponding author}
\footnotetext[2]{\tt The research of the first author was in part supported by a grant from IPM (No. 89050046)}


\begin{abstract}
In this paper we prove that any distance-balanced graph $G$ with $\Delta(G)\geq |V(G)|-3$ is regular.
Also we define notion of distance-balanced closure of a graph and we find distance-balanced closures of trees $T$ with
$\Delta(T)\geq |V(T)|-3$.

\vspace{3mm}
 \end{abstract}

\section{Introduction}

Let $G$ be a graph with vertex set $V(G)$ and edge set $E(G)$.
We denote $|V(G)|$ by $n$.
The set of neighbors of a vertex $v\in V(G)$
is denoted by $N_G(v)$, and $N_G[v]=N_G(v) \cup \{v\}$.
The degree of a vertex $v$ is denoted by $\deg_G(v)$ and minimum degree and maximum degree of $G$
denoted by $\delta(G)$ and $\Delta(G)$, respectively.
The distance $d_G(u,v)$ between vertices $u$ and $v$ is the length of a shortest path
between $u$ and $v$ in $G$. The diameter ${\rm diam}(G)$ of graph $G$ is defined as
$\max \{d_G(u,v): u,v \in V(G)\}$.

For an edge $xy$ of a graph $G$ $W_{xy}^G$ is the set of vertices closer to $x$ than $y$, more formally
$$W_{xy}^G=\{u\in V(G)| d_G(u,x)< d_G(u,y)\}.$$
Moreover, $_xW_{y}^G$ is the set of vertices of $G$ that have equal distances to $x$ and $y$ or
$$_xW_{y}^G=\{u\in V(G)| d_G(u,x)= d_G(u,y)\}.$$
These sets play important roles in metric graph theory, see for
instance \cite{bandel, chepo, djok, epp}. The set
$W_{xy}^G$ also appears in chemical graph theory to
compute the well-investigated Szeged index of a graph (see \cite{grao, gut}), defined as:
$$Sz(G)= \sum_{xy\in E(G)} |W_{xy}^G|. |W_{yx}^G|.$$
{\it Distance-balanced} graphs are introduced in \cite{dbJKR} as graphs for which
$|W^G_{xy}|=|W^G_{yx}|$ for every pair of adjacent vertices $x,y \in V(G)$.

In \cite{dbJKR}, the parameter $b(G)$ of a graph $G$ is introduced
as the smallest number of the edges which can be added to $G$ such
that the obtained graph is distance-balanced. Since the complete
graph is distance-balanced, this parameter is well-defined. We call graph $G$
a {\it distance-balanced closure} of $H$ if $H$ is a spanning
subgraph of $G$ and $|E(G)|= b(H)+|E(H)|$; in other words, a distance-balanced closure of $H$ is a graph $G$ which contains $H$ as a spanning subgraph and has minimum number of edges. As mentioned in \cite{dbJKR}, the  computation of $b(G)$ is quite hard in general but it might be interesting in some special cases. In this paper we compute
 $b(G)$ for all trees $T$ with $\Delta(T)\geq |V(T)|-3$. In Section 2, 
 we compute that distance-balanced closure
 of graphs $G$ with $\Delta(G)= n-1$. In Section 3, and Section 4, we concern graphs $G$ with $\Delta(G)=n-2$
 and $\Delta(G)=n-3$, respectively. Then we compute $b(T)$ for all trees $T$ with $\Delta(T)=n-2$ and $\Delta(T)=n-3$.

Here we mention some more definitions and notations about trees.
Let $P_n$ denoted the path with $n$ vertices. A tree which has
exactly one vertex of degree greater than two is said to be {\it
starlike}. The vertex of maximum degree is called the {\it central
vertex}. We denote by $S(n_1, n_2, \ldots , n_k)$ a starlike tree in
which removing the central vertex leaves disjoint paths $P_{n_1},
P_{n_2}, \ldots , P_{n_k}$. We say that $S(n_1, n_2, \ldots , n_k)$
has branches of length $n_1, n_2, \ldots , n_k$. It is obvious that
$S(n_1, n_2, \ldots , n_k)$ has $n_1 + n_2 + \ldots +n_k + 1$
vertices.
For simplicity a starlike with $\alpha_i$ branches of length $n_i$ ($1\leq i\leq k$) is denoted by
$S(n_1^{\alpha_1}, n_2^{\alpha_2}, \ldots , n_k^{\alpha_k})$.

\section{Distance-balanced graphs with maximum degree $n-1$}

In this section we prove for any graph $G$ with $\Delta(G)=n-1$, the only distance-balanced closure
of $G$ is the complete graph $K_n$.
The following result is very useful in this paper. It is in fact a slight modification of Corollary 2.3 of \cite{dbJKR}.

\begin{thm}\label{dim2} Let $G$ be a graph with diameter at most $2$ and $H$ be a distance-balanced graph
 such that $G$ is a spanning subgraph of $H$. Then $H$ is a regular graph. Moreover,
 every regular graph with diameter at most $2$ is distance-balanced.
\end{thm}

\begin{cor}\label{star}
For every integer $n\geq 1$, the graph $K_{1,m}$ has a unique distance-balanced closure which
is isomorphic to $K_{m+1}$, hence, $b(K_{1,m})= {m+1 \choose 2} -m$.
\end{cor}

\begin{proof} {Let $G$ be a distance-balanced closure of $K_{1,m}$. By Theorem \ref{dim2}, $G$ is a regular
graph and since $K_{1,m}$ has a vertex of degree $m$, $G$ should be
$m$-regular, hence $G \cong K_{m+1}$.}
\end{proof}

The following is an immediate conclusion of Theorem \ref{dim2}.

\begin{cor} Let $G$ be a graph with $n$ vertices and $\Delta (G)=n-1$. Then the graph $G$ has a unique closure which is isomorphic to $K_{n}$.
\end{cor}

\section{Distance-balanced graphs with maximum degree $n-2$}

In this section, we prove that any distance-balanced graph $G$ with $\Delta(G)=n-2$ is a regular graph,
using this we construct a distance-balanced closure of $G$ where $G$ is a tree with this property and then compute $b(G)$.

The following Lemma will be used occasionally in this paper and the proof is easily deduced from the definition of $W^G_{xy}$.
\begin{lem}\label{wn}
{Let $x$ and $y$ be two vertices of a graph $G$, then $W_{xy}^G \cap N_G(y)=\emptyset$. Furthermore,
$N_G(y) \setminus W_{yx}^G \subseteq N_G(x)$.}
\end{lem}
\begin{thm}\label{s2reg}
Let $G=S(2, 1^{m-1})$ be a starlike tree and $H$ be a distance-balanced graph containing $G$ as a spanning subgraph. Then ${\rm diam}(H)\leq 2$, hence, $H$ is an $r$-regular graph for some $m\leq r \leq m+1$.
\end{thm}

\begin{proof}
{ Suppose that the vertices of $G$ are labeled as shown in Figure 1.

\vspace*{.1cm}
{\small $$\epsfbox{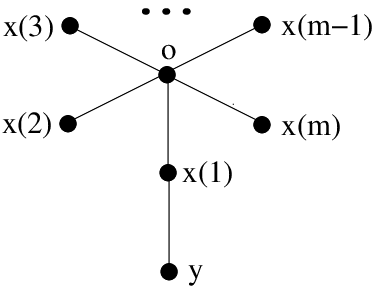}$$}
\begin{center}
{\rm Figure 1}
\end{center}
\vspace*{.1cm}

If $oy \in E(H)$, then $H$ contains $K_{1,m+1}$ as a spanning subgraph, so by Theorem
 \ref{dim2} and Corollary \ref{star}, $H \cong K_{m+2}$.

 So, we may assume that $oy\not \in E(H)$ (and consequently ${\rm diam}(H)\neq 1$).
We prove ${\rm diam}(H)=2$; For this, it is enough to show that $d(y,x_i)\leq 2$ for $i=1,\cdots,m$.
Let $i$ be an integer with $1\leq i\leq n$; From $|W_{x_i o}^H| = |W_{o x_i }^H |$ and Lemma \ref{wn}, we conclude
 that $W_{x_i o}^H \subseteq \{y\}$; We may suppose $W_{x_i o}^H = \emptyset$ (Otherwise, $W_{x_i o}^H = \{y\}$ leads to $x_iy \in E(H)$ and $d_{H}(x_i,y)=1$).
 By Lemma \ref{wn}, this leads to
$$N_H(x_i)=\{o, x_1, x_2, \ldots, x_m\}\setminus \{x_i\},$$
from which, we conclude that $x_1\in N_H(x_i) \cap N_H(y)$ and $d_H(x_i, y)=2$. Hence, ${\rm diam}(H)=2$, as required. The result is now concluded from Theorem \ref{dim2}.
}
\end{proof}

\begin{thm}\label{s2}
Let $G=S(2, 1^{m-1})$ be a starlike of order $m+2$. Then
$$b(G)= \left\{
  \begin{array}{ll}
    \frac{m^2}{2}-1 & \hbox{if $m$ is even;} \\
    {m+1} \choose {2} & \hbox{otherwise.}
  \end{array}
\right.$$
\end{thm}

\begin{proof}
{Let the vertices of $G$ be labeled as in Figure 1, and ${\overline G}$ be a distance-balanced closure of $G$. First,
suppose that $m$ is an odd integer; Since there is no $m$-regular
graph of order $m+2$, by Theorem \ref{s2reg}, ${\overline
G} \cong K_{m+2}$.

Now, suppose that $m$ is an even integer. Let $H=K_{m+2}$ be a complete graph with vertex set
$V(G)$ and $M$ be a complete matching of $H$ which contains the edge $oy$. Then $H\setminus M$,
is an $m$-regular graph with diameter $2$. Hence, by Theorem \ref{dim2} and Theorem \ref{s2reg},
${\overline G}=H\setminus M$, is a distance-balanced closure of $G$ and $b(G)=\frac{m^2}{2}-1$. }
\end{proof}

\begin{cor}
Let $G$ be a connected graph of order $n$, with $\Delta (G)= n-2$ and $H$ be a distance-balanced graph which contains $G$ as a spanning subgraph. Then $H$ is either an $(n-2)$-regular graph or the complete graph $K_n$.
\end{cor}

\begin{proof}
{ In this case $S(2, 1^{n-3})$ is an spanning subgraph of $G$. So, by Theorem \ref{s2reg}, $H$ is either
an $(n-2)$-regular graph or the complete graph $K_{n}$.}
\end{proof}

\section{Distance-balanced graphs with maximum degree $n-3$}

In this section we will prove that every distance-balanced graph with $\Delta(G)=n-3$ is regular.
Moreover, by constructing distance-balanced closure of trees with $\Delta(T)=n-3$ we compute $b(T)$ for these
trees.

\begin{thm}\label{s22reg}
 Let $G=S(2^2 , 1^{m-2})$ be a starlike of order $m+3$ and $H$ be a distance-balanced graph
 which contains $G$ as a spanning subgraph. Then ${\rm diam}(H)\leq 2$, hence, $H$ is an
 $r$-regular graph with $r\geq m$.
\end{thm}

\begin{proof}{ Suppose the vertices of $G$ are labeled as in Figure 2.

\vspace*{.1cm}
{\small $$\epsfbox{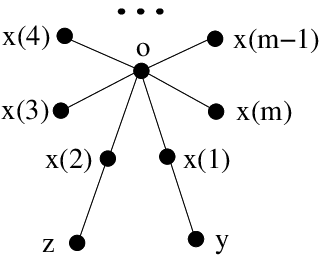}$$}
\begin{center}
{\rm Figure 2}
\end{center}
\vspace*{.1cm}

If either $o y$ or $o z$ be an edge of $H$, then by Theorem \ref{s2reg}, $H$ is an $r$-regular graph with $r\geq m+1$, which proves this theorem.
So, suppose that $oy,oz \not \in E(H)$.

By Lemma \ref{wn}, for each $1 \leq i\leq m$, $W_{x_i o}^H \subseteq \{y, z\}$.
Now, we prove that $\deg_H(x_i)=m$, for each $i=1,\cdots,m$; For this, we consider
 three possible cases:

\begin{description}
\item [Case 1. $|W_{x_i o}^H| = 0$.] Then $W_{x_i o}^H =\emptyset$ and
$W_{o x_i}^H =\emptyset$.
Hence, by Lemma \ref{wn}, $N_H(x_i)=\{o, x_1, x_2, \ldots, x_m\} \setminus \{x_i\}$ and $\deg_H(x_i)=m$.

\item [Case 2. $|W_{x_i o}^H| = 1$.] Without loss of generality we can assume that
 $W_{x_i o}^H =\{y\}$. Then there is an integer $1 \leq j \leq m$ such that
  $W_{o x_i}^H=\{x_j\}$. Since $d_H(o, y)=2$, $yx_i\in E(H)$.
 Hence, using Lemma \ref{wn}, $N_H(x_i)=\{o, y, x_1, x_2, \ldots, x_m\} \setminus \{x_j\}$ and $\deg_H(x_i)=m$.

\item[ Case 3. $|W_{x_i o}^H| = 2$.] We have  $W_{x_i o}^H =\{y, z\}$.
  Since  $d_H(o, y)=d_H(o, z)=2$, we conclude $yx_i, zx_i\in E(H)$. Since $|W_{ox_i}^H| = 2$,
there are integers $j$ and $k$ such that
$W_{o x_i}^H =\{x_j, x_k\}$.

   Hence, by Lemma \ref{wn}, we have $N_H(x_i)=\{o, y, z, x_1, x_2, x_3, \ldots, x_m\} \setminus \{x_i, x_j, x_k\}$ and $\deg_H(x_i)=m$.
\end{description}

Next, we prove that $\deg_H(y)\geq m-3$ and $\deg_H(z)\geq m-3$.
From $\deg_H(x_1)=m$ and Lemma \ref{wn}, it concludes that $|W_{yx_1}^H|\leq 2$, hence, $|W_{x_1y}^H|\leq 2$,
 which means that there are at most two elements in
$N_H(x_1)\setminus N_H[y]$. Using this and Lemma \ref{wn}, we provide
$\deg_H(y)\geq m-3$. With a similar argument, the inequality $\deg_H(z)\geq m-3$ is concluded.

Now, by using $\deg_H(y),\deg_H(z)\geq m-3$, $\deg_H(x_i)\geq m,\,$ $(i=1,\cdots,m)$, and
$oy,oz\not \in E(H)$, we conclude that every two nonadjacent vertices have a common neighbor, provided that $m\geq7$. This means that ${\rm diam}(H)=2$, which proves the result in case $m\geq 7$, using Theorem \ref{dim2}.
 For cases $3\leq m \leq 6$, through a case by case inspection (by using $\deg_H(x_i)\geq m,\,$ $i=1,\cdots,m$,) the same result is obtained.
}\end{proof}

\begin{thm}\label{s22}
For the starlike tree $G=S(2^2, 1^{m-2})$ of order $m+3$, $b(G)= \frac{m^2+m-4}{2}$.
\end{thm}

\begin{proof}{
Let the vertices of $G$ be labeled as in Figure 2 and
${\overline G}$ be a distance-balanced closure of $G$. Now, we are going to construct ${\overline G}$.
Let $H=K_{m+3}$ be a complete graph with the same vertex set as $H$.
Omit the edges of cycles $C_1= x_1 x_2 x_3\ldots x_m x_1$ and
$C_2=oyzo$ from $H$ to obtain ${\overline G}=H\setminus (C_1 \cup C_2)$. Now, ${\overline G}$
 is an $m$-regular graph with diameter 2, which contains $G$ as a spanning subgraph, so by Theorem \ref{s22reg} and Theorem \ref{dim2},
${\overline G}$ is a distance-balanced closure of $G$ and $b(G)= \frac{m^2+m-4}{2}$.
}\end{proof}

\begin{thm}\label{tn3}
Let $G$ be the graph of Figure 3 and $H$ be a distance-balanced graph which contains
$G$ as a spanning subgraph. Then ${\rm diam}(H)\leq 2$, hence, $H$ is a regular graph. Moreover,
$b(G)= \frac{m^2+m-4}{2}$.
\vspace*{.1cm}
{\small $$\epsfbox{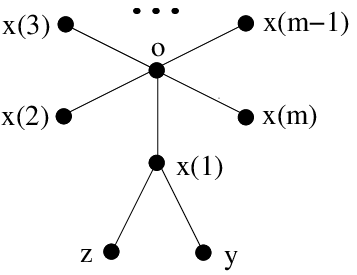}$$}
\begin{center}
{\rm Figure 3}
\end{center}
\vspace*{.1cm}
\end{thm}

\begin{proof}{
If either $oy\in E(H)$ or $oz\in E(H)$, then by Theorem \ref{s2reg}, ${\rm diam}(H)\leq 2$ and $H$ is a regular graph.
So, suppose that neither $oy$ nor $oz$ is in $E(H)$.
Since $|W_{x_iy}^H|=|W_{yx_i}^H|$ and $o\in W_{x_iy}^H$, there exists a vertex $x_i$, $i\neq 1$, such that
$yx_i\in E(H)$. Therefore, graph $H$ contains graph $S(2^2, 1^{m-2})$ as a spanning subgraph and using
Theorem \ref{s22reg}, ${\rm diam}(H)\leq 2$ and $H$ is a regular graph. Furthermore, the graph introduced
in the proof of Theorem \ref{s22}, is also distance-balanced closure of $G$. Hence $b(G)= \frac{m^2+m-4}{2}$.
}
\end{proof}

\begin{thm}\label{s3reg}
Consider the starlike tree $G=S(3, 1^{m-1})$ of order $m+3$
  and let $H$ be a distance-balanced graph which contains $G$ as a spanning subgraph. Then $H$ is an $r$-regular graph for some $m\leq r\leq m+2$.
\end{thm}

\begin{proof}{ Let the vertices of $G$ be labeled as in Figure 4.
\vspace*{.3cm}
{\small $$\epsfbox{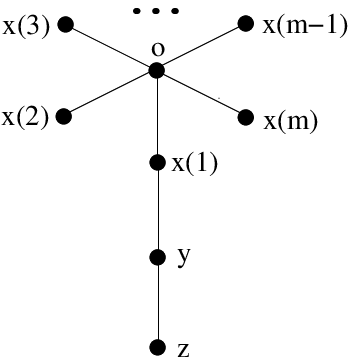}$$}
\begin{center}
{\rm Figure 4}
\end{center}
\vspace*{.1cm}

If either $oy\in E(H)$ or $oz\in E(H)$, then by Theorem \ref{s2reg}, ${\rm diam}(H)\leq 2$ and
$H$ is a regular graph. Suppose that $oy\not\in E(H)$ and $oz\not\in E(H)$.
Since $|W_{yz}^H|=|W_{zy}^H|$, the vertex $z$
is adjacent with at least one vertex in
$\{x_1, x_2, \ldots, x_m\}$. Hence, $H$ contains either the graph of Figure 3 or the graph
$S(2^2, 1^{m-2})$, as a spanning subgraph. So, by Theorem \ref{s22reg} and Theorem \ref{tn3}, ${\rm diam}H\leq 2$
and $H$ is a regular graph, as desired.
}\end{proof}

\begin{cor}
Let $G$ be a connected graph of order $n$ with $\Delta(G)= n-3$. Then every distance-balanced graph $H$ which contains $G$ as a spanning subgraph, is regular.
\end{cor}

\begin{proof}{
Since $\Delta(G)= n-3$, $G$ contains at least one of the graphs $S(2^2, 1^{n-2})$, $S(3, 1^{n-1})$ or the graph shown in Figure 3, as a spanning subgraph. Hence, the result follows from Theorem \ref{s22}, Theorem \ref{s3bg} and
Theorem \ref{tn3}.
}
\end{proof}

\begin{thm}\label{s3bg}
For the starlike tree $G=S(3, 1^{m-1})$ of order $m+3$,
$b(G)=\frac{m^2 + m - 4}{2}$.
\end{thm}
\begin{proof}{
Let the vertices of $G$ be labeled as in Figure 4 and let
${\overline G}$ be a distance-balanced closure of $G$. Now, we are going to construct ${\overline G}$.
Let $H=K_{m+3}$ be a complete graph with the same vertex set as $G$.
Omit the edges of cycles $C_1= x_3 x_4\ldots x_m x_3$ and
$C_2=oyx_2x_1zo$ from $H$ to obtain ${\overline G}=H\setminus (C_1 \cup C_2)$. Then the graph ${\overline G}$ is an $n$-regular graph with diameter 2, which contains $G$ as a spanning subgraph. So by Theorem \ref{s3reg},
${\overline G}$ is a distance-balanced closure of $G$ and $b(G)= \frac{m^2+m-4}{2}$.
}
\end{proof}

{\bf Conclusion.}
In previous sections, we have proved that any connected distance-balanced graph $G$ with $\Delta(G)\geq |V(G)|-3$ is a regular graph, moreover, distanced-closure of such a graph $G$ is a smallest regular graph which contains $G$.
This helped us to find a distance-balanced closure of trees $T$ with $\Delta(T)\geq |V(T)|-3$
and to compute $b(T)$ for such trees.

\end{document}